\documentclass[twoside,english]{elsarticle}
\usepackage[fleqn]{amsmath}
\usepackage[T1]{fontenc}
\usepackage[utf8]{inputenc}
\usepackage{mathtools}
\usepackage{amsthm}
\usepackage{amssymb}
\usepackage{tikz}
\usepackage{pgfplots}
\pgfplotsset{compat=1.13}
\usepackage{todonotes}

\makeatletter
\theoremstyle{plain}
\newtheorem{thm}{\protect\theoremname}
\theoremstyle{remark}
\newtheorem{rem}{\protect\remarkname}

\journal{a journal}

\makeatother

\usepackage{babel}
\providecommand{\remarkname}{Remark}
\providecommand{\theoremname}{Theorem}

\begin{document}

\begin{frontmatter}{}

\title{An analysis of over-relaxation in kinetic approximation\tnoteref{t1}}

\tnotetext[t1]{This work has been supported by ANR EXAMAG project. }

\author[rvt]{Florence Drui\corref{cor1}}

\ead{florence.drui@math.unistra.fr}

\author[rvt]{Emmanuel Franck}

\author[rvt]{Philippe Helluy}

\author[rvt]{Laurent Navoret}

\cortext[cor1]{Corresponding author}

\address[rvt]{Université de Strasbourg \& Inria Tonus}
\begin{abstract}
The over-relaxation approach is an alternative to the Jin-Xin
relaxation method (\citet{jin1995relaxation}) in order to apply the equilibrium
source term
in a more precise way (\citet{coulette2016palindromic,coulette2017palindromic}).
This is also a key ingredient of the Lattice-Boltzmann method for
achieving second order accuracy (\citet{dellar2013interpretation}).
In this work we provide an \hyphenation{ana-ly-sis}analysis of the
over-relaxation kinetic scheme. We compute its equivalent equation,
which is particularly useful for devising stable boundary conditions
for the hidden kinetic variables.
\end{abstract}


\end{frontmatter}{}

\section{Introduction}

\global\long\def\vf{\mathbf{{f}}}
\global\long\def\vg{\mathbf{{g}}}
\global\long\def\vh{\mathbf{{h}}}
\global\long\def\vs{\mathbf{{s}}}
\global\long\def\vV{\mathbf{V}}
\global\long\def\vx{\mathbf{x}}
\global\long\def\v#1{\mathbf{#1}}
\global\long\def\vq{\mathbf{q}}
\global\long\def\vw{\mathbf{w}}
\global\long\def\ddim{D}
\global\long\def\dorder{d}
\global\long\def\normal{\mathbf{n}}
\global\long\def\flux{\mathbf{q}}
\global\long\def\Vu{\mathbf{u}}
\global\long\def\Vvi{\mathbf{v}_{i}}
\global\long\def\wini{\mathbf{v}}
\global\long\def\vu{\v u}

In this note, we are interested in the numerical resolution of the
following system of conservation laws
\begin{equation}
\partial_{t}\vu+\partial_{x}\vf(\vu)=0,\label{eq:conslaw}
\end{equation}
where the unknown is the vector of conservative variables $\vu(x,t)\in\mathbb{R}^{m}$,
depending on a space variable $x$ and a time variable
$t\geq0$. In the first part of the paper we consider the case with no
boundaries ($x \in\mathbb{R}$). In Section \ref{sec:bc}, we will discuss the 
case with boundaries ($x \in [0,1]$). The conservative variables satisfy an initial condition
\begin{equation}
\vu(x,0)=\wini(x).\label{eq:cond_ini}
\end{equation}

The flux $\vf$ is a non-linear function of $\vu$. The system of
conservation laws (\ref{eq:conslaw}) is assumed to be hyperbolic:
for any vector of conservative variables $\vu$, the jacobian of the
flux 
\begin{equation}
\v A(\vu)=\vf'(\vu)\label{eq:defjacob}
\end{equation}
is diagonalizable with real eigenvalues.  

Jin and Xin (\citet{jin1995relaxation}) have proposed an approximation
of (\ref{eq:conslaw}) of the following form
\begin{eqnarray}
\partial_{t}\vw_{\varepsilon}+\partial_{x}\v z_{\varepsilon} & = & 0,\label{eq:relax1}\\
\partial_{t}\v z_{\varepsilon}+\lambda^{2}\partial_{x}\vw_{\varepsilon} & = & \frac{1}{\varepsilon}(\vf(\vw_{\varepsilon})-\v z_{\varepsilon}),\label{eq:relax2}
\end{eqnarray}
with the initial condition
\begin{equation}
\vw_{\varepsilon}(x,0)=\wini(x),\quad\v z_{\varepsilon}(x,0)=\vf(\wini(x)).\label{eq:cond_ini_relax}
\end{equation}
In this formulation $\varepsilon$ is a small positive parameter and
$\lambda$ a constant positive velocity. If $\lambda$ is large enough
(``subcharacteristic'' condition), it can be proved
that $\vw_{\varepsilon}$ tends to $\vu$, the entropy solution of (\ref{eq:conslaw})-(\ref{eq:cond_ini}),
when $\varepsilon$ tends to zero (\citet{jin1995relaxation}).

The advantage of the Jin-Xin formulation is that the partial differential
equations are now linear with constant coefficients and the non-linearity
is concentrated in a simple source term. 

A simple way to solve numerically
the Jin-Xin system is to use a time-splitting algorithm. For advancing
by one time step of size $\Delta t$, one first solves
\begin{eqnarray}
\partial_{t}\vw+\partial_{x}\v z & = & 0,\label{eq:relax1-1}\\
\partial_{t}\v z+\lambda^{2}\partial_{x}\vw & = & 0,\label{eq:relax2-1}
\end{eqnarray}
for a duration of $\Delta t$ (free transport step). Then, 
for the same duration, one solves
 the system of differential equations (relaxation step)
\begin{eqnarray}
\partial_{t}\vw & = & 0,\label{eq:relax1-2}\\
\partial_{t}\v z & = & \frac{1}{\varepsilon}(\vf(\vw)-\v z).\label{eq:relax2-2}
\end{eqnarray}
Both sub-steps admit a simple explicit solution.
Indeed, if we define the functional translation operator $\tau(\Delta t)$ by
\[
(\tau(\Delta t)\v v)(x)=\v v(x-\lambda\Delta t),
\]
then the solution of the free transport step (\ref{eq:relax1-1})-(\ref{eq:relax2-1})
is given by
\[
\left(\begin{array}{c}
\vw(\cdot,t+\Delta t)\\
\v z(\cdot,t+\Delta t)
\end{array}\right)=T(\Delta t)\left(\begin{array}{c}
\vw(\cdot,t)\\
\v z(\cdot,t)
\end{array}\right),
\]
with
\begin{equation}
T(\Delta t)\coloneqq\frac{1}{2}\left(\begin{array}{cc}
\tau(\Delta t)+\tau(-\Delta t) & (\tau(\Delta t)-\tau(-\Delta t))/\lambda\\
\lambda(\tau(\Delta t)-\tau(-\Delta t)) & \tau(\Delta t)+\tau(-\Delta t)
\end{array}\right).\label{eq:free_transport}
\end{equation}
This can be easily obtained noting that the characteristic quantities $\v w / 2 \pm \v z / 2\lambda$ are transported at velocities $\pm \lambda$. 
The solution of the relaxation step is given by
\[
\left(\begin{array}{c}
\vw(\cdot,t+\Delta t)\\
\v z(\cdot,t+\Delta t)
\end{array}\right)=P_{\varepsilon}(\Delta t)\left(\begin{array}{c}
\vw(\cdot,t)\\
\v z(\cdot,t)
\end{array}\right),
\]
with
\begin{equation}
P_{\varepsilon}(\Delta t)\left(\begin{array}{c}
\vw\\
\v z
\end{array}\right)\coloneqq\left(\begin{array}{c}
\vw\\
\vf(\vw)
\end{array}\right)+\exp(-\Delta t/\varepsilon)\left(\begin{array}{c}
0\\
\v z-\vf(\vw)
\end{array}\right).\label{eq:exact_relaxation}
\end{equation}
These operators being defined, we obtain a first-order-in-time approximation of the solution
of (\ref{eq:relax1})-(\ref{eq:relax2})-(\ref{eq:cond_ini_relax})
\[
\left(\begin{array}{c}
\vw_{\varepsilon}(\cdot,\Delta t)\\
\v z_{\varepsilon}(\cdot,\Delta t)
\end{array}\right)=S_{1}(\Delta t)\left(\begin{array}{c}
\v v\\
\vf(\v v)
\end{array}\right)+O(\Delta t^{2}),
\]
with
\begin{equation}
S_{1}(\Delta t)=P_{\varepsilon}(\Delta t)T(\Delta t).\label{eq:lie_scheme}
\end{equation}
The splitting error is of order $O(\Delta t^{2})$, but when this approximation is accumulated on $t/\Delta t$ time steps,
$S_{1}$ is indeed a first order scheme. The Jin-Xin scheme is very robust and can handle shock solutions.
However, for smooth solutions, its  accuracy is not sufficient.

For achieving second order accuracy for smooth solutions, a simple idea would be to
replace the splitting (\ref{eq:lie_scheme}) by a Strang procedure.
We observe that $T(0)=I$, where $I$ is the identity operator. For
$\varepsilon>0$ fixed, we also have $P_{\varepsilon}(0)=I$. However, when
$\varepsilon$ tends to zero, the relaxation step becomes
\[
P_{0}(\Delta t)\left(\begin{array}{c}
\vw\\
\v z
\end{array}\right)=\left(\begin{array}{c}
\vw\\
\vf(\vw)
\end{array}\right)
\]
and we observe that the limit relaxation operator does not satisfy
$P_{0}(0)=I$ anymore. It has become a projection operator
$$P_{0}(0)P_{0}(0)=P_{0}(0). $$
The fact that $P_{0}(0)\neq I$ is the main reason why a Strang splitting
procedure like
\[
S(\Delta t)=T(\frac{\Delta t}{2})P_{\varepsilon}(\Delta t)T(\frac{\Delta t}{2})
\]
 would not lead to a second order scheme in the case $\varepsilon = 0$ (\citet{coulette2016palindromic,coulette2017palindromic}).
The objectives of this paper are:
\begin{enumerate}
\item Recall how to construct a splitting that remains second order when
$\varepsilon=0$.
\item Compute the formal equivalent equation of the resulting scheme.
\item From this equivalent system of partial differential equations construct compatible boundary conditions
ensuring stability and high order.
\item Test the whole approach for a simple hyperbolic problem solved with
a Lattice-Boltzmann Method. 
\end{enumerate}

\section{Over-relaxation scheme}


For constructing a second-order-in-time over-relaxation scheme, a possibility is to perform
a Padé approximation when $\Delta t\simeq0$ of the exponential operator
\[
\exp(-\frac{\Delta t}{\varepsilon})\simeq\frac{1-\frac{\Delta t}{2\varepsilon}}{1+\frac{\Delta t}{2\varepsilon}}=\frac{2\varepsilon-\Delta t}{2\varepsilon+\Delta t},
\]
and to replace the exact relaxation (\ref{eq:exact_relaxation}) step by 
\begin{equation}
R_{\varepsilon}(\Delta t)\left(\begin{array}{c}
\vw\\
\v z
\end{array}\right)\coloneqq\left(\begin{array}{c}
\vw\\
\vf(\vw)
\end{array}\right)+\frac{2\varepsilon-\Delta t}{2\varepsilon+\Delta t}\left(\begin{array}{c}
0\\
\v z-\vf(\vw)
\end{array}\right).\label{eq:exact_relaxation-1}
\end{equation}
We would obtain the same formula by applying a Crank-Nicolson
scheme for approximating the differential equation (\ref{eq:relax1-2})-(\ref{eq:relax2-2}).
Now, we observe that
\begin{equation}
R_{0}(\Delta t)\left(\begin{array}{c}
\vw\\
\v z
\end{array}\right)=\left(\begin{array}{c}
\vw\\
2\vf(\vw)-\v z
\end{array}\right).\label{eq:zero_relax}
\end{equation}
This operator does not depend on $\Delta t$ anymore. We also observe that, as for the
usual Strang splitting procedure, $R_{0}(0)\neq I$.
In addition, $R_{0}$ is no more a projection, but an involutory operator
\[
R_{0}R_{0}=I.
\]
With this observation in mind, we propose the following over-relaxation
scheme $S_{2}(\Delta t)$ for approximating the solution of (\ref{eq:conslaw})-(\ref{eq:cond_ini}).
It is defined by
\begin{equation}
S_{2}(\Delta t)\coloneqq T(\frac{\Delta t}{4}) \, R_{0}\, T(\frac{\Delta t}{2})\, R_{0}\, T(\frac{\Delta t}{4}).\label{eq:scheme_def}
\end{equation}
With this definition, we can check that the over-relaxation scheme
is time-symmetric:
\[
S_{2}(-\Delta t)=S_{2}(\Delta t)^{-1},\quad S_{2}(0)=I.
\]
This property ensures that the over-relaxation scheme 
is second order in time (\citet{hairer2006geometric,mclachlan2002splitting}).
For one single time step we thus have
\[
\left(\begin{array}{c}
\vu(\cdot,\Delta t)\\
\vf(\vu(\cdot,\Delta t))
\end{array}\right)=S_{2}(\Delta t)\left(\begin{array}{c}
\v v\\
\vf(\v v)
\end{array}\right)+O(\Delta t^{3}),
\]
where $\vu$ is the exact solution of (\ref{eq:conslaw})-(\ref{eq:cond_ini}).

\section{Equivalent equation}

In this section, we will compute the equivalent equation of the over-relaxation
scheme. The objective is to derive a system of partial differential
equations satisfied by the approximations of $\vw$ and $\v z$ when
$\Delta t$ tends to zero. Of course, if $\v z = \vf(\vw)$ at the initial time, we expect $\vw$ and $\v z$
to satisfy
\[
\partial_{t}\vw+\partial_{x}\vf(\vw)=O(\Delta t^{2}),\quad\v z-\vf(\vw)=O(\Delta t^{2}).
\]

A more interesting question is to find a partial differential equation
satisfied by the approximation of the flux $\v z.$ This is important
in practice in order to construct
 stable boundary conditions to be applied
to $\v z$, or for designing schemes that remain second order at the
boundaries.

Let $\vw$ and $\v z$ denote now the numerical solution given by
the second order scheme (\ref{eq:scheme_def}). We have
\begin{multline}
\frac{1}{\Delta t}\left(\begin{array}{c}
\vw(\cdot,t+\frac{\Delta t}{2})-\vw(\cdot,t-\frac{\Delta t}{2})\\
\v z(\cdot,t+\frac{\Delta t}{2}))-\v z(\cdot,t-\frac{\Delta t}{2}))
\end{array}\right) = \\
\frac{1}{\Delta t}(S_{2}(\frac{\Delta t}{2})-S_{2}(-\frac{\Delta t}{2}))\left(\begin{array}{c}
\vw(\cdot,t)\\
\v z(\cdot,t)
\end{array}\right),\label{eq:def_s2}
\end{multline}
where the over-relaxation scheme $S_{2}$ has an explicit form given
by (\ref{eq:free_transport}), (\ref{eq:zero_relax}) and (\ref{eq:scheme_def}).
We can perform a Taylor expansion of both sides of \eqref{eq:def_s2} when $\Delta t$
tends to zero. By symmetry considerations, the first order terms vanish.
An essential point is that $S_{2}(0)=I$. After simple but long calculations,
we obtain the following result:
\begin{thm}
Let $\vw$ and $\v z$ be smooth solutions of the time marching algorithm
(\ref{eq:def_s2}). Let us define the flux error $\v y$ by
\[
\v{y\coloneqq z-f(w)}.
\]
Then, up to second order terms in $\Delta t$, $\vw$ and $\v y$
are solutions of the following (non conservative) hyperbolic system
of conservation laws:
\begin{equation}
\partial_{t}\left(\begin{array}{c}
\vw\\
\v y
\end{array}\right)+\left(\begin{array}{cc}
\vf'(\vw) & 0\\
0 & -\vf'(\vw)
\end{array}\right)\partial_{x}\left(\begin{array}{c}
\vw\\
\v y
\end{array}\right)=0.\label{eq:limit_wz}
\end{equation}
\end{thm}
\begin{rem}
The equivalent equation (\ref{eq:limit_wz}) shows that the over-relaxation
scheme tends to propagate the conservative variables $\vw$ and the
flux error $\v y$ with opposite wave velocities. This gives hints to
build stable boundary conditions on $\v z$. Roughly speaking, at
an inflow boundary for $\vw$, one should impose $\vw$ and not $\v y$,
while at an outflow boundary for $\vw$ one should impose $\v y$
and not $\v w$. Numerical experiments that confirm this heuristic
are given in the next section.
\end{rem}
\begin{rem}
	Generally, when the relaxation operator is a projection, the equivalent equation is only available for $\vw$. See for instance
	\cite{dubois2008equivalent,sportisse2000analysis}
\end{rem}
\begin{rem}
The kinetic speed $\lambda$ does not appear in the equivalent equation
(\ref{eq:limit_wz}). It appears in the $O(\Delta t^{2})$ terms,
which are complicated. We do not know yet how to perform the stability
analysis of these terms. In practice, if $\lambda$ is too small, the scheme becomes unstable. This indicates that a subcharacteristic condition
still has  to be satisfied.
\end{rem}

\section{Numerical method}

\subsection{Transport model}

In this section, we describe a numerical discretization of $S_{2}$
in the simple case where $m=1$, $\vu=u\in\mathbb{R}$ and $\vf(\vu)=f(u)=cu$,
$c>0$. We thus solve a simple transport equation at velocity $c>0$
\[
\partial_{t}u+c\, \partial_{x}u=0.
\]
We assume that $x\in[0,1].$ We also provide an initial condition
and a boundary condition at the left point
\[
u(x,0)=v(x),\quad u(0,t)=v(-ct).
\]
where $v:\mathbb{R}\rightarrow\mathbb{R}$ is a given function. It can be checked that the exact solution of this initial-boundary
value problem is
\[
u(x,t)=v(x-ct).
\]
With this transport equation, we can associate its over-relaxation system
with the approximated conservative data $\vw=w\in\mathbb{R}$ and
the approximated flux $\v z=z\in\mathbb{R}.$

\subsection{Numerical discretization}

For the numerical discretization, we consider a positive integer $N$
and define the space step and grid points by
\[
\Delta x=\frac{1}{N+1},\quad x_{i}=i\Delta x,\quad i=0\ldots N+1.
\]
The grid points $i=0$ and $i=N+1$ are the border points of the interval
$[0,1]$, where the boundary conditions are applied. We
consider an approximation of $w$ and $z$ at the grid points $x_{i}$
and times $t_{n}=n\Delta t$
\[
w_{i}^{n}\simeq w(x_{i},t_{n}),\quad z_{i}^{n}\simeq z(x_{i},t_{n}).
\]
The initial data are exactly sampled at the grid points
\[
w_{i}^{0}=v(x_{i}),\quad z_{i}^{0}=cv(x_{i}).
\]
Like in the Lattice Boltzmann Method (\citet{chen1998lattice}) we
choose a special time step
\[
\Delta t=\frac{4\Delta x}{\lambda}.
\]
This choice ensures that the transport operator $T(\Delta t/4)$ only
involves exact shift operators. For instance, the translation operator
is approximated here by
\[
(\tau(\Delta t/4)w)(x_{i},t_{n})\simeq w_{i-1}^{n}.
\]
Thus, the transport step reads:
\begin{equation}
\begin{aligned}
	w_i^{n+1/4} & =  \frac{w_{i-1}^{n} + w_{i+1}^{n}}{2} + \frac{z_{i-1}^n - z_{i+1}^n}{2\lambda}, \\
	z_i^{n+1/4} & =  \frac{z_{i-1}^{n} + z_{i+1}^{n}}{2} + \lambda \frac{w_{i-1}^n - w_{i+1}^n}{2}. 
\end{aligned}
\label{eq:discrete_scheme}
\end{equation}
Note that, these discrete equations are equivalent to:
\begin{align}
	z_i^{n+1/4} - \lambda w_i^{n+1/4} & =  z_{i+1}^{n} - \lambda w_{i+1}^{n}, \label{eq:discrete_scheme_f-}\\
	z_i^{n+1/4} + \lambda w_i^{n+1/4} & =  z_{i-1}^{n} + \lambda w_{i-1}^{n}. \label{eq:discrete_scheme_f+}
\end{align}
In practice, we also use the fact that $T(\Delta t/2)=T(\Delta t/4)T(\Delta t/4)$.

\subsection{Boundary conditions, relaxation \label{sec:bc}}

\subsubsection{Boundary conditions\label{subsec:Boundary-conditions}}

Let us assume that at the beginning of a time step, for instance at time $t_n$, we know $w_{i}^{n}$
and $z_{i}^{n}$ for $i=0\ldots N+1$. The transport operator $T(\Delta t/4)$
can be applied to internal grid points $x_{i}$, corresponding to
indices $i=1\ldots N.$ Thus, using \eqref{eq:discrete_scheme},  it is possible to compute $w_{i}^{n+1/4}$,
$z_{i}^{n+1/4}$ for $i=1\ldots N$. At the left boundary $i=0$,
one information is missing for computing $w_{0}^{n+1/4}$ and $z_{0}^{n+1/4}.$
According to the previous analysis, it is natural to impose the boundary
condition on $w$ (because it is an inflow boundary) at the middle
of the time step (in order to respect the time-symmetry):
\[
w(0,t_{n}+\frac{\Delta t}{8})=v(-c(t_{n}+\frac{\Delta t}{8})).
\]
It is discretized by
\begin{equation}
\frac{w_{0}^{n}+w_{0}^{n+1/4}}{2}=v(-c(t_{n}+\frac{\Delta t}{8})), \label{eq:w-left-bounds}
\end{equation}
which provides the missing relation and enables to compute $w_{0}^{n+1/4}$. Then, from \eqref{eq:discrete_scheme_f-}, the value of $z_{0}^{n+1/4}$ can be computed:
		\[
			z_0^{n+1/4} - \lambda w_0^{n+1/4} = z_1^n - \lambda w_1^n.
		\]

At the right boundary, we will test several approaches: an ``exact''
strategy, a ``Dirichlet'' strategy on $y=z-cw$ or a ``Neumann''
strategy.

\paragraph*{Exact strategy}

Since we know the analytical solution, we can impose the values
of $w_{N+1}$ 
 given by the exact solution. Of course
this method cannot be generalized to more complex equations and solutions.
In addition, we expect it to generate oscillations, because
the boundary condition is not compatible with an outflow boundary.
As for the left boundary, we write
		\[
			\frac{w_{N+1}^{n}+w_{N+1}^{n+1/4}}{2}=v(1-c(t_{n}+\frac{\Delta t}{8})).
		\]

\paragraph*{Dirichlet strategy on $y$}

In this method, we simply apply the condition $y=0$ at the middle
of the time step. This gives
\[
\frac{z_{N+1}^{n}+z_{N+1}^{n+1/4}}{2}-c\frac{w_{N+1}^{n}+w_{N+1}^{n+1/4}}{2}=0,
\]
which provides the missing relation.
		For instance, using this relation and expression \eqref{eq:discrete_scheme_f+}, 
		one obtains
		\[
			w_{N+1}^{n+1/4} =  \frac{1}{\lambda + c} \left( \lambda w_{N}^{n} -c w_{N+1}^{n}\right) 
			+ \frac{1}{\lambda + c} \left( z_{N}^{n} + z_{N+1}^{n}\right).
		\]

\paragraph*{Neumann strategy on $y$}

The last method consists in imposing the condition $\partial_{x}y(L,t)=0$
at the right boundary. Formally, up to second order, this is equivalent
to impose $\partial_{t}y(L,t)=0$ or $y(L,t)=y(L,0)=0.$ The missing
relation is obtained from
\[
z_{N+1}^{n+1/4}-cw_{N+1}^{n+1/4}=z_{N}^{n+1/4}-cw_{N}^{n+1/4}.
\]
		Once again, using this relation and expression \eqref{eq:discrete_scheme_f-}, 
		 one now has the following relation for $w_{N+1}^{n+1/4}$
		\begin{multline*}
			w_{N+1}^{n+1/4} = \frac{1}{2(\lambda+c)} \left[ \left( 2 \lambda w_{N}^{n} + (\lambda +c) w_{N+1}^{n} - (\lambda-c) w_{N-1}^{n}\right) 
			+ \right. \\
                   \left. \left( 2 z_{N}^{n} -\frac{\lambda + c}{\lambda} z_{N+1}^{n} -\frac{\lambda - c}{\lambda} z_{N-1}^{n}\right) \right].
		\end{multline*}

\subsubsection{Relaxation}

The relaxation operation, as stated above, consists in replacing in
each cell $(\vw,\v z)$ by $(\vw,2\vf(\vw)-\v z)$. We emphasize
that the relaxation is also performed in the boundary cells $i=0$
and $i=N+1.$

\subsection{Numerical results}

We test the above scheme and boundary conditions with $x\in[0,1]$,
$t\in[0,t_{\max}]$ and the following exact solution
\[
u(x,t)=\exp(A(x-\alpha-ct)^{2}).
\]
We may also  impose a (non-physical) flux disequilibrium $y=z-f(w)\neq 0$ at the initial time.  
The initial value of $y$ is given at time $t=0$ by
\[
	y =  B\exp(A(x-\beta+ct)^{2}).
\]

We test the three boundary approaches proposed in Section \ref{subsec:Boundary-conditions}.
We check the stability and order of the scheme. The error is measured
by the discrete $L^{2}$ norm
\[
e_{\Delta x}^{n}=\sqrt{\Delta x\sum_{i=0}^{N+1}\big(w_{i}^{n}-u(x_{i}-cn\Delta t)\big)^{2}+\big(z_{i}^{n}-cu(x_{i}-cn\Delta t)\big)^{2}}.
\]

\begin{figure}
	\begin{center}
		\includegraphics[width=0.9\textwidth]{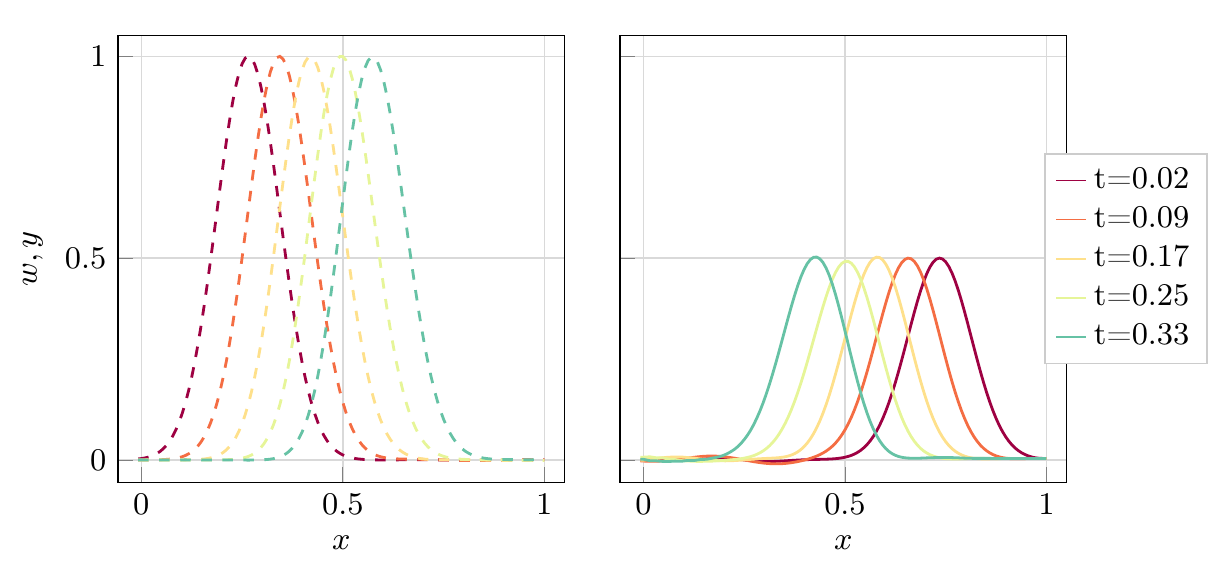}
		\caption{Transport of the $w$ (dashed lines) and $y = z-f(w)$ (plain lines) quantities.}
		\label{fig:y-transport-illustr}
	\end{center}
\end{figure}

\begin{figure}
	\begin{center}
		\includegraphics[width=0.9\textwidth]{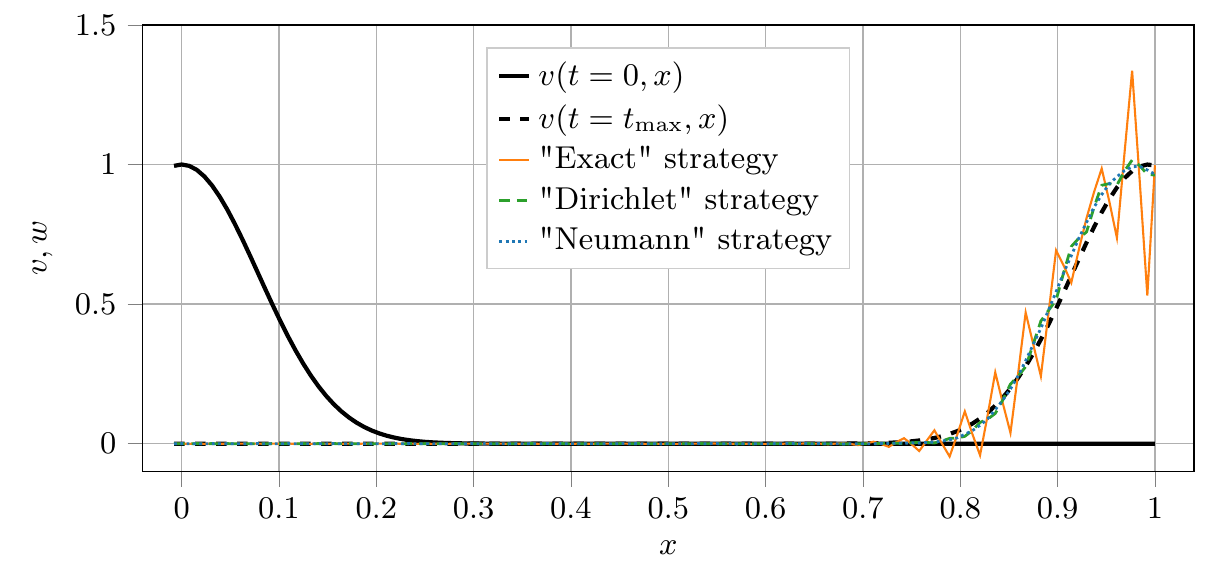}
		\caption{Initial state and comparison of the final states for the  transport 
		equation with Gaussian initial profile, $\Delta x = 2^{-7}$.}
		\label{fig:illustr-transport-test}
	\end{center}
\end{figure}
\begin{figure}
	\begin{center}
		\includegraphics[width=0.9\textwidth]{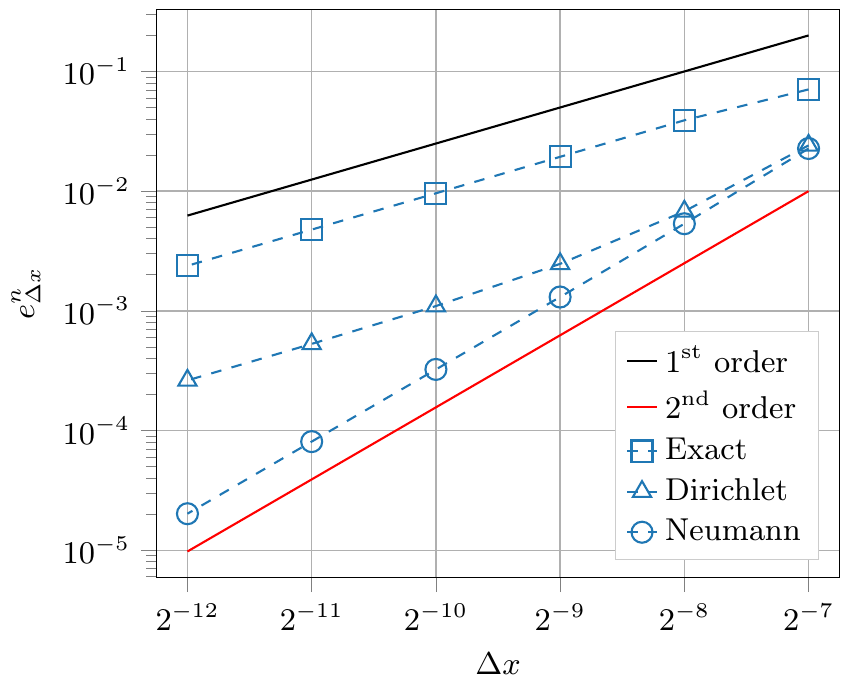}
		\caption{Convergence study for the transport equation with Gaussian initial 
		profile. Comparison of Exact, Dirichlet and Neumann strategies.}
		\label{fig:convergence-transport-test}
	\end{center}
\end{figure}

In Figure \ref{fig:y-transport-illustr}, we first show an illustration of 
the propagation of the quantity $y = z-f(w)$
with the following numerical parameters:
\[
c=1,\quad 
\lambda=2,\quad 
t_{\max}=0.33,\quad
\alpha=0.25,\quad 
\beta = 0.75\quad 
A=80,\quad
\text{and}\quad 
B=1/2.
\]
With this choice and sufficiently small final time, the boundary condition has no influence. 
One checks numerically that $w$ propagates with velocity $u$,  while $y$ propagates
with velocity $-u$, which agrees with the equivalent equation \eqref{eq:limit_wz}.

An illustration of the numerical results using the three strategies for boundary conditions is given in Figure \ref{fig:illustr-transport-test}. For this test, we have imposed the following numerical parameters:
\[
c=1,\quad 
\lambda=2,\quad 
t_{\max}=1,\quad
\alpha=0,\quad 
\beta = 0\quad 
A=80,\quad
\text{and}\quad 
B=0.
\]
One can see that the \textit{Exact} strategy generates oscillations at the right boundary. 
The \textit{Dirichlet} strategy generates weaker oscillations that are not amplified with time.
The \textit{Neumann} strategy does not generate any oscillation.

The convergence results are shown in Figure \ref{fig:convergence-transport-test}. One can see that the
 \textit{Exact} strategy and the  \textit{Dirichlet} strategy are first order accurate, while the  
 \textit{Neumann} strategy is second order accurate.
 The best choice for ensuring stability and second order accuracy seems to be
 the \textit{Neumann} strategy.
 
 \begin{rem}
 	With the  \textit{Dirichlet} strategy, we impose that $y = z-cw = 0$ at
	the boundary, while this equality may not be satisfied exactly inside the 
	domain. This creates small discontinuities that may explain the loss of accuracy. 
 \end{rem}

\section{Conclusion}

In this short note, we have derived the equivalent equation of the over-relaxation
kinetic scheme. The equivalent equation reveals that the conservative
variable and the flux error propagate in opposite directions. This
allows us to determine natural boundary conditions for the over-relaxation
scheme. Numerical experiments confirm the stability and accuracy of
these boundary conditions. In a forthcoming work, we will extend the
approach to more complex non-linear systems and to higher dimensions.
It also important to incorporate in the over-relaxation method a dissipative 
mechanism in order to compute discontinuous solutions without oscillations.
\section{Bibliography}

\bibliographystyle{elsarticle-num-names}
\addcontentsline{toc}{section}{\refname}\bibliography{crmeca}

\end{document}